\documentclass[12pt]{amsart}
\usepackage{amssymb,amsmath}

\theoremstyle{plain}
\newtheorem{theorem}{Theorem}[section]

\newtheorem{proposition}[theorem]{Proposition}

\theoremstyle{definition}
\newtheorem{definition}[theorem]{Definition}

\theoremstyle{remark}

\begin{document}
\title[Wavelet filter functions]
{Wavelet filter functions, the matrix completion problem, 
and projective modules over $C(\mathbb T^n)$}
\thanks{The first author was supported by a research grant 
from the National University of Singapore. The research of the
second author was
supported in part by National Science Foundation grant DMS99-70509.}
\author[J.A.~Packer]{Judith A. Packer}
\address{Department of Mathematics, National University of Singapore, 
10 Kent Ridge Crescent, Singapore 119260, Republic of Singapore}
\address{Department of Mathematics, University of Colorado at Boulder,
Boulder, CO 80309-0395, U.S.A.}
\email{packer@euclid.colorado.edu}
\author[M.~A.~Rieffel]{Marc A.~Rieffel}
\address{Department of Mathematics, University of California at 
Berkeley, Berkeley, CA 94720, U.S.A.}
\email{rieffel@math.berkeley.edu}
\subjclass
%[2000]
{Primary 46L99; Secondary 42C40, 46H25}
%%%%%%%%%%%%%%%%%%%%%%%%%%%%%%%%%%%%%%%%%%%%%%%%%%%%%%%%%%%%%%%%%%%%%%%%%%%%%
\begin{abstract}
We discuss how one can use certain filters from signal processing to describe 
isomorphisms between  certain projective
$C(\mathbb T^n)$-modules. Conversely, we show how cancellation 
properties for finitely generated projective modules over $C(\mathbb 
T^n)$
can often be used to prove the existence of continuous high pass filters,
of the kind needed for multi-wavelets, 
corresponding to a given continuous low-pass filter. However, we also give an 
example of a continuous low-pass filter for which it is impossible 
to find corresponding continuous high-pass filters. In this way we 
give another approach to the solution of the matrix completion 
problem for filters of the kind arising in wavelet theory.
\end{abstract}
%%%%%%%%%%%%%%%%%%%%%%%%%%%%%%%%%%%%%%%%%%%%%%%%%%%%%%%%%%%%%%%%%%%%%%%%%%%%
\maketitle

A key technique in wavelet theory is to use
suitable low-pass filters to construct scaling functions and their
multi-resolution analyses, and then to use corresponding high-pass filters to
construct the corresponding wavelets.  
Thus once one has used a low-pass filter
to construct a scaling function, it is important to be able to find associated
high-pass filters.  (They are not unique.)  It has been pointed out several
times in the literature that, given a continuous low-pass filter for dilation
by $q>2$, there is in general no continuous selection function for the
construction of associated continuous high-pass filters, because of the
existence of topological obstructions.  (See the substantial discussion in
section 2 of \cite{BrJ}, and the references cited there.)  It is known
\cite{My}, \cite{Wo}, \cite{BrJ} that associated high-pass filters which are
{\bf measurable} will always exist.  (See also the general Hilbert-space
treatment of the existence of wavelets given in \cite{BC}.)  The problem of
finding appropriate high-pass filters given a low-pass filter is closely
related to the more general problem of 
completing an $n\times n$ unitary matrix
whose entries take values in $C(\mathbb T^n)$ when given its first row.  Thus
this problem is sometimes called the ``matrix completion problem'' for
wavelets.

The main theorem in the literature concerning the existence of high-pass
filters in the multivariate case was first given by Gr\"ochenig \cite{Gr} for
dyadic dilation matrices.  It is given a nice exposition in section 3.4 of
\cite{My2}.  The generalization to general dilation matrices is given an
attractive treatment in theorem 5.15 of \cite{Wo}.  This theorem states the
following.  Let $A$ be an $n\times n$ dilation matrix with
$|\text{det}(A)|\;=\;q$.  Let $m_0$ be a low-pass filter function on $\mathbb
T^n$ associated with $A$ which is what 
is called $r$-regular.  If $2q-1\;>\;n$,
then there exist associated high-pass filters 
which are $r$-regular.  (The case
in which $q = 2$ and $n \geq 3$ is not covered by this theorem, but it is
well-known that the case $q=2$ can always be treated by simpler methods, which
we indicate in Section 3 below.)  The methods of proof are analytical.  The
condition of $r$-regularity implies that the corresponding filters (i.e.
transfer functions) are infinitely differentiable.  
In Section 3 of this paper
we show, using the methods of cancellation for projective modules (or vector
bundles), that if we are given a {\bf continuous} low-pass filter function
$m_0$ on $\mathbb T^n$ corresponding to $A$, and if $2q-1\;>\;n$, then there
will always exist {\bf continuous} high-pass filters
$m_1,\;m_2,\;\cdots,m_{q-1}$ associated to $m_0$.

For most applications one does want the filters to have a certain degree of
smoothness.  Our emphasis on only requiring that the filters be continuous is
mainly to make clear that the issues discussed here are basically topological.
Smoothness can always be restored, as discussed at the end of the proof of our
main theorem.

We do not give an explicit algorithm for constructing the high-pass filters,
which is not surprising since the proof 
of the cancellation theorem quoted from
\cite{Hu} is not algorithmic in nature, and we already know from the
topological obstructions mentioned above that it is impossible to find an
algorithm into which one can plug low-pass filters $m_0$ and read out in a
continuous way corresponding continuous high pass filters.  In both our
situation and that of theorem 5.15 of \cite{Wo} and 
the earlier references, the
issue of constructing, for general dilation matrices, low-pass filters which
generate continuous scaling functions is left open.

Finally in the fourth section of the paper we give a fairly explicit example,
for $n=5$ and $q = 3$, of a continuous low-pass filter for which it is {\bf
impossible} to construct a family of corresponding continuous high-pass
filters.  (Thus the condition $2q-1\;>\;n$ is best possible for $n=5$.)  Our
construction uses classical facts about vector bundles over spheres.  We note
that it is shown in \cite{GR} that if one is content with wavelet frames, and
with the number of wavelets exceeding $q$, then one can always find compactly
supported wavelets of arbitrarily high smoothness.

The second author pointed out in a talk given in 1997 \cite{Rie2} that, at
least for filters whose transfer functions are continuous, one can use
``projective modules'' over $C(\mathbb T^n)$ to gives an attractive framework
for wavelet theory.  In the present paper we make more explicit the
relationship between filters and projective modules over $C(\mathbb T^n)$.  It
is essential for our results that our filters have {\bf continuous} transfer
functions.  (We will refer to such filters simply as ``continuous filters''.)
In a subsequent paper \cite{pacrief} we will directly use non-free projective
modules to give a new variation on the construction of wavelets.

In the first two sections of the paper we establish our basic framework, and
reformulate well-known facts to show that from a collection of continuous
filters for an $n\times n$ integer dilation matrix $A$, one can construct an
explicit isomorphism between a certain module over $C(\mathbb T^n)$ determined
by $A$ and the free modules $\oplus_{i=1}^q [C(\mathbb T^n)]$, where
$q\;=\;|\text{det}(A)|.$

We remark that we will make extensive use of ``inner products'' with values in
algebras such as $C(\mathbb T^n)$.  Such inner products have been widely used
in some related parts of harmonic analysis.  See \cite{La} and references
therein.  These inner products have also been used to some extent within the
signal processing literature, 
often under the name of ``bracket product''.  See
for example the discussion of bracket products for Weyl-Heisenberg frames in
chapter 11 of \cite{Ca}.  These inner products are used intensively in the
first sections of \cite{Rie5}, which can appropriately be viewed as a
discussion of windowed Fourier transforms 
for lattices in ${\mathbb R}^n$ which
are skew to the standard basis.  These inner products seem to have been used in
only a few places in the wavelet literature.  See, for example, the bracket
product defined in equation 2.6 of \cite{RS}, 
and associated references.  These
inner products could profitably be used more widely, 
since many of the formulas
in wavelet theory have attractive reformulations in terms of them, as we will
see in part $2$ below.  These inner products will play an important role in
\cite{pacrief}.

The most rapidly accessible introduction to 
multivariate wavelets with which we
are familiar is that given in \cite{Str}, and we will have it in mind during
our discussion below.  Another attractive exposition is that given in section
5.2 of \cite{Wo}.

It seems likely that our techniques can be applied in the case of
multi-wavelets, as described for example in \cite{AK}, but we have not explored
this possibility.

The first author would like to thank Professors Lawrence Baggett and Dana
Williams for many useful discussions on the topics discussed in this paper, and
for their great hospitality towards her and her family during her visit to the
University of Colorado at Boulder and Dartmouth College, respectively, during
her sabbatical year.

%%%%%%%%%%%%%%%%%%%%%%%%%%%%%%%% 
\section{Actions of finite subgroups on
$C(\mathbb T^n),$ free modules, and module frames}
%%%%%%%%%%%%%%%%%%%%%%%%%%%%%%%%%%%%%%%%

Let $A$ be an $n \times n$ integer matrix, and set $q = |\text{det}(A)|$.
Assume that $q \geq 2$.  (Eventually we will require that $A$ is a dilation
matrix.)  Let $\Gamma = A\mathbb Z^n$.  It is easily seen that $\Gamma$ is a
subgroup of $G = \mathbb Z^n$ of finite index $q$.  Any subgroup of $\mathbb
Z^n$ of finite index arises in this way.  For the next few paragraphs we can
forget about $A$ and just keep $\Gamma$.

Let ${\mathcal B} = C_f(G)$ denote the convolution algebra of complex-valued
functions on $G$ of finite support.  (This is where finite impulse response
(FIR) filters reside.)  Similarly we set ${\mathcal A} = C_f(\Gamma)$.  In the
evident way ${\mathcal A}$ is a subalgebra of ${\mathcal B}$, and so we can
view ${\mathcal B}$ as a module over ${\mathcal A}$ using convolution.  Let
$\{p_1, \cdots, p_q\}$ be a set of coset representatives for the cosets of
$\Gamma$ in $G$.  For each $j$ let $e_j$ denote the delta-function at $p_j$.
Each element, $f$, of ${\mathcal B}$ is the sum of its restrictions to the
cosets, and so is uniquely of the form $f = \sum h_j\ast e_j$ for $h_j$'s in
${\mathcal A}$.  This means, by definition, that ${\mathcal B}$ is a {\bf free}
${\mathcal A}$-module of rank $q$, with the $e_j$'s serving as a module basis.
The mapping $f \mapsto (h_j) \in {\mathcal A}^q$ (column vectors with entries
in ${\mathcal A}$) defines an ${\mathcal A}$-module isomorphism from ${\mathcal
B}$ to ${\mathcal A}^q$ .

In this kind of situation there is a natural ${\mathcal A}$-valued ``inner
product'' (or ``bracket product'' as mentioned in the introduction) on
${\mathcal B}$, whose use goes back at least to corollary 4.7 of \cite{Rie1a}.
We recall the natural involution on ${\mathcal B}$ defined by $f^*(m) = {\bar
f}(-m)$.  For $f, g \in {\mathcal B}$ we form the convolution $f \ast (g^*)$
and then restrict it to $\Gamma$ (i.e.  by ``down-sampling'' or
``decimating'').  That is, we set $$ \langle f, g \rangle_{\mathcal A}(\gamma)
= (f \ast (g^*))(\gamma) = \sum_G f(m){\bar g}(m-\gamma) , 
$$ where $\gamma \in
\Gamma$.  (This inner product is ${\mathcal A}$-linear in its first variable.)
It is easily checked that $\{e_j\}$ is an ``orthonormal'' basis for this inner
product.  All of this structure can be extended to various completions of
${\mathcal A}$ and ${\mathcal B}$, 
such as those for the $\ell^1$-norm, and the
C*-norm.  See \cite{La} and references therein for the general theory of such
inner products.

We will work mostly with the Fourier-transformed version of the above picture.
Each of the dual groups $\hat G$ and $\hat \Gamma$ is isomorphic to ${\mathbb
T}^n$, but we must distinguish carefully between these two dual groups, and it
is only $\hat G$ which we will identify with ${\mathbb T}^n$.  We view
${\mathbb T}^n$ as ${\mathbb R}^n/{\mathbb Z}^n$, and as in \cite{Str} we use
the exponentials $e^{2\pi i n\cdot t}$ in defining Fourier series.  Let $F$ be
the subgroup of ${\hat G} = {\mathbb T}^n$ consisting of the characters of $G$
which take value $1$ on all of $\Gamma$.  Then $F$ can be viewed as the dual
group of $G/\Gamma$, and it is a finite group of order $q$.  Its elements,
viewed as characters of $G$, act by pointwise multiplication to give
automorphisms of the convolution algebra ${\mathcal B}$.  It is easily seen
that ${\mathcal A}$ consists exactly of the elements of ${\mathcal B}$ which
are left fixed by all the automorphisms from $F$.  As the Fourier transform of
${\mathcal B}$ we take the completion $C({\mathbb T}^n)$, consisting of the
complex-valued continuous functions on ${\mathbb T}^n$, with pointwise
multiplication, and adjoint given by pointwise complex conjugation.  We will
denote it again by ${\mathcal B}$.  
(This should not cause confusion.)  In this
picture the action of the group $F$ on ${\mathcal B}$ consists simply of
translation by the elements of $F$.  Then (the completion of) ${\mathcal A}$
becomes just the subalgebra of ${\mathcal B}$ consisting of the functions in
$C({\mathbb T}^n)$ which are invariant under translation by the elements of
$F$.  We denote it again by ${\mathcal A}$.  In this picture the ${\mathcal
A}$-valued inner product on ${\mathcal B}$ is then easily seen to be given by
$$ \langle f, g \rangle_{\mathcal A}(x) = q^{-1}\sum_{w \in F} f(x-w){\bar
g}(x-w).  $$ That the factor $q^{-1}$ gives the correct normalization can be
checked by considering the identity element of ${\mathcal B}$.

Of course, in the Fourier picture ${\mathcal B}$ 
is still free as an ${\mathcal
A}$-module.  To see this directly, consider the dual group, $\hat F$ ($=
G/\Gamma)$, of $F$.  Extend each element of $\hat F$ to a character of $\hat
G$.  (This corresponds to choosing coset representatives for $\Gamma$ in $G$;
there is no canonical choice of extensions in general).  Let $\{e_j\}$ denote
this collection of characters.  A simple calculation shows that $\{e_j\}$ is
orthonormal for the ${\mathcal A}$-valued inner product, and that $$ f = \sum
\langle f, e_j \rangle_{\mathcal A} e_j $$ for every $f \in {\mathcal B}$.
Thus $\{e_j\}$ is an ${\mathcal A}$-module basis for ${\mathcal B}$.  It is
very natural to define an ${\mathcal A}$-valued inner product on ${\mathcal
A}^q$ by $$ \langle (a_j), 
(c_j)\rangle_{\mathcal A} = \sum a_j{\bar c}_j .  $$
The mapping $f \mapsto (\langle f, \;e_j\rangle_{\mathcal A})$ 
is ``isometric''
for this inner product.  Of course, we should expect that there will be many
other orthonormal ${\mathcal A}$-module bases for ${\mathcal B}$, and this will
be the subject of the next section.

We summarize some of the above discussion with:

\begin{proposition} \label{thm free} Let $F$ be a finite subgroup of ${\mathbb
T}^n$ of order $q$, and let ${\mathcal A}$ denote the subalgebra of ${\mathcal
B} = C({\mathbb T}^n)$ consisting of the functions which are invariant under
translation by the elements of $F$.  View ${\mathcal B}$ as a module over
${\mathcal A}$.  Then ${\mathcal B}$ is a free module over ${\mathcal A}$ with
$q$ generators.  \end{proposition}

In Section 4 we will need to deal with ${\mathcal A}$-modules which are not
free, but are (finitely generated) projective.  These will also be of central
importance for a paper about wavelets presently under preparation
\cite{pacrief}.  By definition, 
a projective module is (isomorphic to) a direct
summand of a free module.  Let us view our free module as ${\mathcal A}^m$ for
some integer $m$, with ``standard basis'' $\{e_j\}$ Then this means that there
is an $m \times m$ matrix $P$ with entries in ${\mathcal A}$ which is a
projection, that is $P^2 = P$, such that our projective module $\Xi$ is of the
form $\Xi = P{\mathcal A}^m$.  In this (C*-algebraic) setting a standard
argument (see 5Bb in \cite{We}) shows that $P$ can be adjusted so that it is
also ``self-adjoint'' in the evident sense.  Then set $\xi_j = Pe_j$ for each
$j$.  For any $\eta \in \Xi$ we have $$ \eta = P\eta = P(\sum \langle \eta, e_j
\rangle_{\mathcal A}e_j) = \sum \langle \eta, e_j \rangle_{\mathcal A} Pe_j \\
= \sum \langle P\eta, e_j \rangle_ {\mathcal A} \xi_j $$ $$ = \sum \langle
\eta, Pe_j \rangle_{\mathcal A}\xi_j = \sum \langle \eta, \xi_j
\rangle_{\mathcal A}\xi_j .  $$ There is no reason to expect that the $\xi_j$'s
will be independent over ${\mathcal A}$, much less orthonormal.  But anyone
familiar with wavelets will feel comfortable about referring to the $\xi_j$'s
as a ``module frame'' for $\Xi$, as the second author did in \cite{Rie2}.  In
the more general setting of projective modules over C*-algebras the
reconstruction formula \begin{equation} \label{eq pmf} \eta\;=\;\sum \langle
\eta, \xi_j \rangle_{\mathcal A}\;\xi_j.  \end{equation} appears, in different
notation, in \cite{Rie}.  This indicates how the $\mathcal A$-valued inner
products are quite useful in discussing frames.  In particular they are useful
in discussing biorthogonal wavelet bases.  The theory of module frames has been
developed extensively in \cite{LaFr98}.

The ${\mathcal A}$-valued inner products we defined earlier are a special cases
of ones associated with ``conditional expectations'', as noted in lemma 1.1 and
proposition 4.17 of \cite{Rie1a}.  Frames for modules connected with
conditional expectations have played an important role in other situations.
See the extensive development by Watatani in \cite{Wa}, and the references
therein.  He uses the terminology ``quasi-bases'' instead of ``frames''.  Our
situation is a special case of proposition 2.8.1 of \cite{Wa}, in which $F$ is
a finite group of cardinality $q$ acting freely on a compact Hausdorff space
$X,$ and the conditional expectation from $C(X)$ onto the fixed point algebra
is shown to be of the kind for which a quasi-basis will exist.  But we remark
that the projective module need not be free in this more general situation.
For example, consider the two-element group acting on the two-sphere by the
antipodal map.  (See situation 2 in \cite{Rie3}.)

%%%%%%%%%%%%%%%%%%%%%%%%%%%%%%%%%%%% 
\section{Filter functions and free $C(\mathbb T^n)$ modules} 
%%%%%%%%%%%%%%%%%%%%%%%%%%%%%%%%%%%%%%

In this section we briefly review the definition of filters, and their
relationships to scaling functions and wavelets.  We then use these filter
functions to construct explicit isomorphisms between certain finitely generated
projective modules over $C({\mathbb T}^n)$ and free $C({\mathbb T}^n)$-modules.

We now review the definition of filter functions corresponding to an arbitrary
integer dilation matrix.  In part we do this because it is important to be
careful with the bookkeeping for the constants involved.  We state everything
in the Fourier transform picture.  We also assume that our filters are
continuous, so that the equations hold everywhere, and not just a.e.  We will
not state in detail the regularity conditions (``Cohen's condition'', etc.)
since we do not explicitly need them.  We refer to \cite{Str} and \cite{Wo} for
precise statements.

\begin{definition} \label{def lpf} Let $A$ be an 
integral dilation matrix, that
is, an integer matrix all of whose eigenvalues have modulus strictly greater
than $1,$ and let $q\;=\;|\text{det}(A)|.$ A continuous function $m_0 \in
{\mathcal B} = C({\mathbb T}^n)$ is called a {\bf low-pass filter} (or
``mask'') for dilation by $A$ if:  \renewcommand{\labelenumi}{(\roman{enumi})}
\begin{enumerate} \item $m_0(0)\;=\;q$, \quad \quad (``low-pass'' condition,
equation 1.17 of \cite{Str}) \item $\langle 
m_0,\; m_0\rangle_{\mathcal A} = q$
\quad \quad (equation 1.16 of \cite{Str}), \item $m_0(x)\;\not=\;0$ for $x$ in
a sufficiently large neighborhood of $0$ well-related to $A$ (``Cohen's
condition", etc.).  \end{enumerate} \end{definition}

It is shown in \cite{Str}, \cite{Wo} that if $m_0$ is a low-pass filter
corresponding to dilation by $A$, then on setting $B = A^T$ and $$
\Phi(x)\;=\;\Pi_{n=1}^{\infty}[q^{-1}m_0(B^{-n}(x))], \quad \quad
\mbox{(equation 1.47 of \cite{Str})} $$ we find that $\Phi$ is the Fourier
transform of a scaling function $\phi\;\in\;L^2(\mathbb R^n)$ corresponding to
dilation by $A.$ In order to guarantee pointwise convergence of the product
defining $\Phi$ and to ensure that $\Phi$ is continuous, we need only make the
mild assumption that $m_0$ is $C^1$ at $0.$ The scaling function and its
corresponding multiresolution analyses are then used to construct a family of
$q-1$ wavelets corresponding to dilation by $A$.  We note that it has recently
been shown by L.  Baggett and K.  Merrill in \cite{BM} that multi-resolution
analyses exist for every $n\times n$ integral dilation matrix $A$; and only
very recently has it been shown by 
M. Bownik and D. Speegle in \cite{BS} that
in the special $2\times 2$ case, to every integer dilation matrix one can
associate a scaling function whose Fourier transform is smooth and compactly
supported.

In what follows we will never need to use condition (iii) of Definition 2.1.
So when we say ``low-pass filter'' below we will mean just conditions (i) and
(ii).  But it is
an open question, for general dilation matrices, as to how often
one can ensure that condition (iii) does hold in such a way that $m_0$ will
yield a good scaling function.

\begin{definition} Let $A$ be an integral dilation matrix with
$|\text{det}(A)|\;=\;q,$ and let $m_0 \in {\mathcal B}$ be a continuous low
pass filter for $A$.  We say that a family $m_1,\;\cdots,\;m_{q-1} \in
{\mathcal B}$ is a {\bf high-pass filter family} corresponding to the low-pass
filter $m_0$, if $$ \langle m_j,\; m_k \rangle_{\mathcal A} = q\delta_{jk}
\quad \quad \mbox{(equation 1.24 of \cite{Str})} $$ for
$0\;\leq\;j,k\;\leq\;q-1$.  \end{definition}

For the reader's reference we remark that the wavelets
$\psi_1,\cdots,\psi_{q-1} \in\;L^2(\mathbb R^n)$ associated to the filter
functions $\{m_0,m_1,\;\cdots,\;m_{q-1}\}$ are given by the formula 
$$
\hat{\psi_i}(x)\;=\;q^{-1}m_i(B^{-1}x)\Phi(B^{-1}x) 
\quad \quad \mbox{(equation
1.39 of \cite{Str})} $$ 
for $1\leq\;i\leq\;q-1$.  Thus if all of the filter
functions are continuous, and if $m_0$ has sufficient regularity so that the
(Fourier transform of the) scaling function $\Phi$ is continuous, as mentioned
in the sentence after the definition of $\Phi,$ then (the Fourier transforms
of) the multiwavelets will be continuous too.

It is clear that if we rescale all of the $m_j$'s by $q^{-1/2}$ then we obtain
an orthonormal family in ${\mathcal B}$.  
But it will be more convenient for us
to rescale the inner product too, in a way which leads to some traditional
formulas.  We define our rescaled inner product by
$\langle\cdot,\;\cdot\rangle'_{\mathcal A} =
q\langle\cdot,\;\cdot\rangle_{\mathcal A}$, so that $$ \langle
f,\;g\rangle'_{\mathcal A}(x) = \sum_{w \in F} f(x - w){\bar g}(x - w) $$ for
$f, \;g \in {\mathcal B}$.  Suppose now that we have a low-pass filter $m_0$
and a corresponding family $m_1,\;\cdots,\;m_{q-1} \in {\mathcal B}$ of
high-pass filters.  We renormalize these by setting $h_j = q^{-1}m_j$ for $0
\leq j \leq q-1$.  This renormalization corresponds exactly to the
renormalization factor in equation 1.35 of \cite{Str}.  It is clear that the
$h_j$'s then form an orthonormal set in ${\mathcal B}$ for our new inner
product $\langle\cdot,\;\cdot\rangle'_{\mathcal A}$, and that $h_0(0) = 1$.  
In
all that follows we now use this new inner product.

We now use a traditional argument to show that such an orthonormal set will
always be an ${\mathcal A}$-module basis for ${\mathcal B}$.

\begin{theorem} \label{thm quasib} Let ${\mathcal A}$ and ${\mathcal B}$ be
defined as above in terms of the finite subgroup $F$ of order $q$.  Let
$h_0,\;\cdots,\;h_{q-1}$ be an orthonormal set in ${\mathcal B}$.  Then it is
an ${\mathcal A}$-basis for ${\mathcal B}$.  \end{theorem}

\begin{proof}  The crux of the matter is to show that the
$h_j$'s ``span'' ${\mathcal B}$.  Label the elements of $F$ by integers $i$
with $0 \leq i \leq q-1$.  
For each $x \in {\mathbb T}^n$ define a $q \times q$
matrix $U(x) = (u_{ji}(x))$ by $u_{ji}(x) = h_j(x-w_i)$.  The orthogonality
condition for the $h_j$'s says that for each $x$ $$ \delta_{jk} = \sum_i
h_j(x-w_i){\bar h}_k(x-w_i) = \sum_i u_{ji}(x){\bar u}_{ki}(x) .  $$ (See
equation 1.41 of \cite{Str}.)  But this says exactly that
$U(x)U(x)^\ast\;=\;I$.  Because we are in finite dimensions, this implies that
$U(x)$ is unitary, so that also $U(x)^\ast U(x)\;=\;I$.  That is, $$ \sum_i
{\bar h}_i(x-w_j)h_i(x-w_k) = \delta_{jk}.  $$ (See equation 1.43 of
\cite{Str}.)  We apply this for $w_k = 0$ to obtain $$ \left(\sum \langle f,\;
h_i\rangle_{\mathcal A}h_i)\right)(x) = \sum_i\left(\sum_j f(x-w_j){\bar
h}_i(x-w_j)\right)h_i(x) = $$ $$ \sum_j f(x-w_j)\sum_j {\bar h}_i(x-w_j)h_i(x)
= f(x).  $$
\end{proof}

\medskip We remark that $h_0$ gives the first row of the unitary matrix $U$
above.  Finding a corresponding set of high-pass filters, so
$h_1,\;\cdots,\;h_{q-1}$, corresponds to finding the other rows of this matrix
(with entries in $C({\mathbb T}^n)$) such that it is unitary.  This is why the
problem of finding high-pass filters is often referred to as the ``matrix
completion problem''.

%%%%%%%%%%%%%%%%%%%%%%%%%%%%%%%%%% 
\section{When low-pass filters have high-pass filter families} 
%%%%%%%%%%%%%%%%%%%%%%%%%%%%%%%%%%%

As discussed in the introduction, 
it is important to find high-pass filters for
a given low-pass filter.  This is not always possible if one requires the
filters to be continuous, as we will show in the next section.  But here we
give some positive results.  
We will use our Proposition \ref{thm free} and the
cancellation theorem for certain finitely generated projective modules over
$C(\mathbb T^n)$ to prove that if $n<2q-1,$ then given any continuous low
pass-filter $m_0$ defined on $\mathbb T^n$ corresponding to any integer
$n\times n$ dilation matrix $A$ with $|\text{det}(A)|\;=q,$ it is possible to
construct not only measurable, but in fact {\bf continuous} high-pass filters
corresponding to the low pass filter $m_0$.  We do not need the condition of
$r$-regularity, which implies that $m_0$ is infinitely differentiable.  (See
the third sentence of the proof of Theorem 5.15 of \cite{Wo}.)

We now state the most general form of our theorem.  \begin{theorem} \label{thm
contshpf} Fix $n\in \mathbb N,$ and let $F$ be a finite subgroup of ${\mathbb
T}^n$ of order $q$.  Let ${\mathcal A}$ and ${\mathcal B}$ be defined as done
earlier, in terms of $F$.  Suppose that $m_0 \in {\mathcal B}$ and that $m_0$
satisfies conditions (i) and (ii) of Definition 2.1.  Suppose that either $q =
2$, or $n<2q-1.$ Then there exists a family of continuous high-pass filters
$m_1,\;m_2,\;\cdots,m_{q-1} \in {\mathcal B}$ corresponding to $m_0.$ In
particular, if $n \leq 4$ then for any $q$ it will be true that any continuous
low-pass filter will have a corresponding family of continuous high-pass
filters.

If a continuous low-pass filter $m_0$ has a corresponding family of continuous
high-pass filters, and if $m_0$ is infinitely differentiable, then $m_0$ has a
corresponding family of high-pass filters which are infinitely differentiable.
\end{theorem}

\begin{proof}  We use the notation of Proposition \ref{thm
free}.  Set $h_0 = q^{-1}m_0$, so that $\langle h_0,\; h_0\rangle'_{\mathcal A}
= 1$.  Viewing ${\mathcal A}$ as an ${\mathcal A}$-module over itself, we
construct an ${\mathcal A}$-module map $\sigma:\;{\mathcal
A}\;\rightarrow\;{\mathcal B}$ by $$ \sigma(a)\;=\;ah_0.  $$ We note that
$\sigma$ preserves the ${\mathcal A}$-valued inner products:  for any $a, b \in
{\mathcal A}$ we have $$ \langle \sigma(a),\; \sigma(b) \rangle'_{\mathcal
A}\;= \;\langle ah_0,\;bh_0 \rangle'_{\mathcal A}\; = a\langle
h_0,\,h_0\rangle_{\mathcal A}{\bar b}\;=\; \langle a,\; b \rangle_{\mathcal A}
.  $$ Hence $\sigma$ gives an ${\mathcal A}$-module injection of ${\mathcal A}$
into ${\mathcal B}$.  Thus $\sigma({\mathcal A})$ is projective, and so it has
a complementary module, $L$, such that $$ {\mathcal B}\;=\;\sigma(\mathcal
A)\;\oplus\;L .  $$ We can take $L$ to be orthogonal to $\sigma(\mathcal A)$
with respect to the inner product described above.  To see that this orthogonal
complement $L$ exists, note that the ``orthogonal'' projection of $\mathcal B$
onto $\sigma(\mathcal A)$ is given by $f \mapsto \langle f,\; h_0
\rangle'_{\mathcal A}h_0$, so that $$ L\;= \;\{f\;-\; \langle f,\; h_0
\rangle'_{\mathcal A}h_0 :  f \in {\mathcal B}\} .  $$

Suppose that a high-pass filter family $m_1,\;m_2,\;\cdots,m_{q-1} \in
{\mathcal B}$ exists for $m_0$, and set $h_j = q^{-1}m_j$ for each $j$.  Then
from the equation of Definition 2.2 we see that the $h_j$'s for $j \geq 1$ will
be an orthonormal family in $L$.  From Theorem 2.3 we can deduce that they will
actually form an orthonormal basis for $L$, so that $L$ must be a free
$\mathcal A$-module.  Thus to show that a high-pass filter family exists, we
need to show that $L$ is a free $\mathcal A$-module, for then we can obtain an
orthonormal basis.  (Use, for example, the proof of proposition 2.1 of
\cite{Rie}.)  We can then multiply by $q$ to obtain the desired $m_j$'s.  Note
that we have $$ {\mathcal A}\;\oplus\; L \;\cong \;{\mathcal
A}\;\oplus\;{\mathcal A}^{q-1}.  $$ (This says that $L$ is ``stably-free''.)
Thus to show that $L$ is free we need to be able to ``cancel'' one copy of
$\mathcal A$, so that $L \cong {\mathcal A}^{q-1}$.

We treat first the case in which $q = 2$.  This has a simple solution, as is
widely seen in the wavelet literature.  In this case the group $F$ has only two
elements.  Let $w$ denote the non-identity element of $F$.  Choose a continuous
function $\tau$ on ${\mathbb T}^n$ such that $|\tau(x)| = 1$ and $\tau(x+w) =
-\tau(x)$ for all $x \in {\mathbb T}^n$.  This can be done, for example, by
choosing $\tau$ to be a character on ${\mathbb T}^n$ such that $\tau(w) = -1$.
Let $h_0$ as above be given.  Define $h_1$ by $h_1(x) = \tau(x){\bar
h}_0(x+w)$.  Then a simple standard calculation shows that the pair $h_0$,
$h_1$ is an orthonormal set in $\mathcal B$.  We can now apply Theorem 2.3 to
conclude that $h_1$ is a basis for $L$.  (We remark that this case is a special
case of the fact that on any compact space any stably-free line bundle is free.
This is because line-bundles are determined by their first Chern class; see
Example 4.55 of \cite{RW}.  But stably equivalent vector bundles have the same
Chern class, by theorem 16 .4.2 of \cite{Hu}.)

To treat the other case we use a theorem of Swan (theorem 1.6.3 of \cite{Ro})
which shows that for any compact space $X$ the projective modules over $C(X)$
correspond to the complex vector bundles over $X$.  One direction of this
correspondence consists of assigning to a vector bundle its $C(X)$-module of
continuous cross-sections.  This enables us to use the facts about cancellation
of vector bundles which are given in \cite{Hu}.  Our $\mathcal A$-valued inner
products correspond to ``Hermitian metrics'' on vector bundles.

Suppose now that $n\;<\;2q-1$.  Then $<\frac{n}{2}>\;\leq\;q-1$, where $<x>$
denotes as in \cite{Hu} the least integer greater than or equal to $x.$ Then
from theorem 8.1.5 of \cite{Hu} we can deduce immediately that we can cancel
$\mathcal A$, so that $L$ is a free module of rank $q-1$.  We have thus proved
the existence of the desired family of continuous high-pass filters associated
to $m_0.$

Suppose now that $m_0$ is infinitely differentiable, and that we have obtained
a corresponding family $m_1,\;m_2,\;\cdots,m_{q-1}$ of continuous high-pass
filters, perhaps by use of the first part of this theorem.  We can uniformly
approximate the $m_j$'s for $j \geq 1$ arbitrarily closely by infinitely
differential functions, say $g_1,\;g_2,\;\cdots,g_{q-1}$.  These functions need
not be orthogonal.  But we can try to apply to them a ``Gram-Schmidt" process
using $\langle \cdot,\; \cdot \rangle_{\mathcal A}$.  The only care needed to
make this work is that the approximations must be close enough so that, if
$f_j$ denotes the orthogonal projection (defined much as in the early part of
the proof of this theorem) of $g_j$ into the orthogonal complement of the span
of the new $m_1,\;m_2,\;\cdots,m_{j-1}$, then $\langle f_j,\; f_j
\rangle_{\mathcal A}$ must still be close enough to $1$ so that $(\langle
f_j,\; f_j \rangle_{\mathcal A})^{-1/2}$ exists and is an infinitely
differentiable function.  For then we can ``normalize" $f_j$ to obtain the new
$m_j$.  
\end{proof}

\medskip On the other hand, if $m_0$ is a FIR filter (i.e.  a trigonometric
polynomial) then it is quite another matter to determine whether corresponding
high-pass filters can be found which are FIR filters.  In this connection see
the discussion of the matrix completion problem in the non-unitary case, but
for Laurant polynomials, given in section 8 of \cite{JM1}.  It employs the
Quillen-Suslin solution of the Serre conjecture.  See also \cite{JM2}.

Of course, all of this discussion is not of much use except for those dilation
matrices for which there exists a continuous low-pass filter which will produce
a scaling function.  As mentioned earlier, it does not seem to be known for
which dilation matrices such a filter always exists, though one can always find
a measurable such filter, \cite{BrJ}.

%%%%%%%%%%%%%%%%%%%%%%%%%%%%%%%%%%%%%%%%%%%%%%%%
\section{A low-pass filter which does not have a full set of high-pass filters}
%%%%%%%%%%%%%%%%%%%%%%%%%%%%%%%%%%%%%%%%%%%%%%%%

We now construct an example of an integral 
dilation matrix $A$ on $\mathbb R^5$
with a given continuous, or even smooth, 
low-pass filter on $\mathbb R^n$ which
does not have a corresponding family of continuous high-pass filters.  We
emphasize that our essential assumption is that the filters which we consider,
as functions on ${\mathbb Z}^n$, decay at infinity fast enough that their
Fourier transforms on ${\mathbb T}^n$ are continuous.  This condition is
satisfied by almost all filters in practical use.

Our example is based on the fact that on ${\mathbb T}^5$ there are complex
vector bundles which are stably free but not free (see below).  This does not
happen on ${\mathbb T}^n$ for $n \le 4$, since cancellation holds as discussed
in the course of the proof of Theorem 3.2.  (Our construction would also work
for $n > 5$.)

Accordingly, we consider ${\mathbb Z}^5$.  To be specific, and to make contact
with wavelet theory, let $A$ be the dilation matrix 
$$ 
A\;=\;\left(
\begin{array}{cc} 0 & 3 \\ I_4 & 0 \end{array} \right), 
$$ 
so that
$\text{det}(A) = 3$.  Let $\Gamma = A({\mathbb Z}^5)$ so that $\Gamma$ has as
 a generating set 
$\{3e_1,e_2,e_3,e_4,e_5\}$, where the $e_j$'s are the standard
generators for ${\mathbb Z}^5$.  We view the convolution $C^*$-algebra
${\mathcal B} = C^*({\mathbb Z}^5)$ as a module over its subalgebra ${\mathcal
A} = C^*(\Gamma)$.  As a module, ${\mathcal B}$ is free, 
of rank~$3$, with (one
possible) module basis given by the $\delta$-functions supported at $0$, $e_1$
and $2e_1$.

We view all of this in the Fourier picture.  Thus we let ${\mathcal B} =
C({\mathbb T}^5)$.  Of course ${\mathcal A} = C^*(\Gamma)$ is also isomorphic
to $C(\mathbb T^5)$ but we must distinguish carefully between it and
$C^*({\mathbb Z}^5)$.  We let $F$ be the subgroup of elements of ${\mathbb
T}^5$ which, viewed as characters on ${\mathbb Z}^5$, have value $1$ on
$\Gamma$.  So $F$ has $3$ elements, and consists of the order-$3$ subgroup of
the first copy of ${\mathbb T}$ in ${\mathbb T}^5$.  Then $C^*(\Gamma)$
corresponds to the subalgebra, ${\mathcal A}$, of ${\mathcal B} = C({\mathbb
T}^5)$ consisting of functions invariant under translation by $F$.  The action
of ${\mathcal A}$ on ${\mathcal B} = C({\mathbb T}^5)$ is by pointwise
multiplication.  We use the ${\mathcal A}$-valued inner product on ${\mathcal
B}$ as before, defined by $$ \langle f,g\rangle'_{\mathcal A}(x) = \sum_{w \in
F} {\bar f}(x-w)g(x-w).  $$ Of course ${\mathcal B}$ is free of rank~$3$ over
${\mathcal A}$, with orthonormal module basis obtained by renormalizing by
$3^{-1/2}$ the set $\{1,{\hat e}_1,({\hat e}_1)^2\}$, 
where ${\hat e}_1$ is the
character of ${\mathbb T}^5$ corresponding to the character from $e_1$ on the
first copy of ${\mathbb T}$ in ${\mathbb T}^5$.

As before, we take our low-pass filter to be renormalized, so given by a
function $h_0 \in {\mathcal B}$ such that $ \langle h_0,h_0\rangle'_{\mathcal
A} = 1$, and $h_0(0) = 1$.  One wants to find corresponding high-pass filters,
$h_1$ and $h_2$, such that $$ \langle h_j,\;h_k\rangle_{\mathcal A} =
\delta_{jk}\hskip 0.5 in j,k = 0,1,2.  $$ Then $\{h_0,h_1,h_2\}$ will be an
``orthonormal'' basis for ${\mathcal B}$ as ${\mathcal A}$-module by Theorem
2.3.

We now show how to construct an $h_0$ for which it is impossible to find
corresponding continuous $h_1$ and $h_2$.  To see what is involved, let $[h_0]$
denote the ${\mathcal A}$-submodule of ${\mathcal B}$ generated by $h_0$.  Then
$[h_0]$ is a free module of rank~$1$.  We let $[h_0]^{\perp}$ denote the
orthogonal complement of $[h_0]$ in ${\mathcal B}$ for the ${\mathcal
A}$-valued inner product.  Note that if $h_1$ and $h_2$ exist, then they form a
module basis for $[h_0]^{\perp}$, and thus $[h_0]^{\perp}$ is a free ${\mathcal
A}$-module of rank~$2$.  Thus to find our desired example, it suffices to find
$h_0$ such that $[h_0]^{\perp}$ is not a free module.  (However $[h_0]^{\perp}
\oplus [h_0] = {\mathcal B}$ will be free, so that $[h_0]^{\perp}$ will be
``stably free'', and so will represent the same element of the $K$-group
$K_0(T^5)$ as the free module of rank~$2$.  Thus the phenomenon we seek cannot
be detected by the $K_0$-group.)

We base our construction of $h_0$ on classical facts about the homotopy groups
of spheres and unitary groups, which show that on the $5$-sphere $S^5$ one can
construct complex vector-bundles having the properties we seek for ${\mathbb
T}^5$.  We view $S^5$ as the standard unit-sphere in ${\mathbb C}^3$.  We will
follow closely some constructions given in section $24.2$ of \cite{Ste} but we
use slightly more concise notation.  Following $24.2$, we view $S^4$ as that
equator of $S^5$ given as $$ S^4 = \{(v,ir):  v \in {\mathbb C}^2,\ r \in
{\mathbb R},\ \|v\|^2 + r^2 = 1\}.  $$ For simplicity we will label these
points just by $(v,r)$.  We view $v$ as a column vector, and let $v^*$ denote
the corresponding row vector with complex conjugate entries.  Then $vv^*$ is a
$2 \times 2$ matrix, self-adjoint, of rank~$1$ (unless $v=0$) with range in the
subspace spanned by $v$.  Let ${\mathcal U}_2$ denote the $2 \times 2$ unitary
group.  We let $U_0$ be the function from $S^4$ to ${\mathcal U}_2$ defined by
$$ U_0(v,r) = I_2 - 2(1+ir)^{-2}vv^*.  $$ This is formula $5$ of $24.2$.
Steenrod provides in $24.3$ of \cite{Ste} a proof that $U_0$ is not
path-connected through unitaries to a constant map from $S^4$ to ${\mathcal
U}_2$.  The idea of the proof is basically that $U_0$ is closely related to the
suspension of the famous Hopf map from $S^3$ to $S^2$ which generates the
homotopy group $\pi_3(S^2) \cong \mathbb Z$.  (Recall that ${\mathcal
S}{\mathcal U}_2$ is homeomorphic to $S^3$.)  Note that $U_0$ can equally well
be considered to be a unitary element of the $C^*$-algebra $M_2(C(S^4))$ of $2
\times 2$ matrices over $C(S^4)$; we will frequently take this point of view.

There is a traditional way, described in section $18.1$ of \cite{Ste}, to use a
map such as $U_0$ to construct a vector bundle over $S^5$.  Because we are
aiming at ${\mathbb T}^5$, we use a simple variant of this construction to
obtain instead a vector bundle over ${\mathbb T} \times S^4$.  This variant is
described in the course of the proof of Theorem $8.4$ of \cite{Rie5}, where,
without additional complication, it is seen to work easily for non-commutative
(unital) $C^*$-algebras.  We identify vector bundles with their modules of
continuous cross-sections, which are projective modules.

For the immediate purposes of this construction it is convenient to view
${\mathbb T}$ as the interval $[0,1]$ with ends identified.  We view
$C({\mathbb T} \times S^4)$ accordingly.  Then the projective module (alias
vector bundle) over $C({\mathbb T} \times S^4)$ determined by $U_0$, denoted
$X(U_0)$, is the vector space of continuous functions $$ X(U_0) = \{F:  [0,1]
\rightarrow (C(S^4))^2:  F(1) = U_0F(0)\}, $$ with the elements of $C({\mathbb
T} \times S^4)$ acting by pointwise multiplication.  We define an inner product
on $X(U_0)$ with values in ${\mathcal C} = C({\mathbb T} \times S^4)$ by $$
\langle F,G\rangle_{\mathcal C}(s) = F(s)^*G(s), $$ where $G(s)$ is viewed as a
column vector, etc.  It follows immediately from lemma $8.10$ of \cite{Rie5}
that if $X(U_0)$ were a free module, then $U_0$ would 
be path-connected through
unitaries to the constant map on $S^4$ with value $I_2$.  
Since this is not the
case, the module $X(U_0)$ is not free.

However, crucial to our purposes is the fact that the direct sum of $X(U_0)$
with the free module of rank~$1$ {\em is} free (of rank~$3$).  From lemma
$8.6$ of \cite{Rie5} the direct sum of $X(U_0)$ with the free module of
rank~$1$ comes by applying the above construction to 
$$ 
\left(
\begin{array}{cc} U_0 & 0 \\ 0 & I_1 \end{array} \right) = U_0 \oplus I_1 
$$
instead of to $U_0$.  From lemma $8.10$ of \cite{Rie5} the fact that the
direct sum is free is equivalent to the fact that $U_0 \oplus I_1$ is
path-connected to $I_3$ through unitaries in $M_3(C(S^4))$.  In order to try
to obtain an explicit formula for our low-pass filter $m_0$, we need an
explicit path.  From the details given in \cite{Ste} it is not difficult to
see how to produce one.  The key is that $U_0 \oplus I_1$ has the following
factorization:  
$$ 
\left( \begin{array}{cc} U_0(v,r) & 0 \\ 0 & 1
\end{array} \right) = \left( \begin{array}{cc} I_2 - (1-ir)^{-1}vv^* & v \\
-b(ir)v^* & ir \end{array} \right)^* \left( \begin{array}{cc} I_2 -
(1+ir)^{-1}vv^* & v \\ c(ir)v^* & ir \end{array} \right), 
$$ 
where $b$ and
$c$ will be defined below.  Each of the two factors is in ${\mathcal U}_3$,
and so the ${}^*$ applied to the first factor could equally well be
${}^{-1}$.  The two factors extend to functions defined on all of $S^5$
except for one point.  Specifically, as in equation~$2$ of $24.2$ of
\cite{Ste}, we set $$ \phi^+(v,\xi) = \left( \begin{array}{cc} I_2 -
(1+{\bar \xi})^{-1}vv^* & v \\ b(\xi)v^* & \xi \end{array} \right) $$ for
$(v,\xi) \in S^5$, with $b(\xi) = (1+\xi)(1+{\bar \xi})^{-1}$.  We note that
$\phi^+$ is not defined at the south pole $(0,-1)$.  We also see that the
first factor above (before taking its inverse) is just $\phi^+(v,ir)$.  In
the same way, following equation~$4$ of $24.2$ of \cite{Ste}, we set $$
\phi^-(v,\xi) = \left( \begin{array}{cc} I_2 - (1-{\bar \xi})^{-1}vv^* & v
\\ c(\xi)v^* & \xi \end{array} \right) $$ for $(v,\xi) \in S^5$, with
$c(\xi) = (1-\xi)(1-{\bar \xi})^{-1}$.  We note that $\phi^-$ is not defined
at the north pole $(0,1)$.  We also see that the second factor above is just
$\phi^-(v,ir)$.  Each of $\phi^+$ and $\phi^-$ has values in ${\mathcal
U}_3$.

We note that $\phi^+(0,1) = I_3$, while $\phi^-(0,-1) = I_2 \oplus (-I_1)$.
This latter matrix is, of course, path-connected to $I_3$.  If we adjust for
this, and move $S^4$ gradually up to the north pole for $\phi^+$, and down to
the south pole for $\phi^-$, we obtain a path through unitaries from $U_0
\oplus I_1$ to $I_3$.  
We can explicitly implement this as follows.  We let $t$
denote the parameter for our path, ranging over $[0,1]$.  For $\phi^+$ we take
the straight-line path $(1-t)ir + t$ from $ir$ to $1$.  
Given $(v,ir)$, we must
scale $v$ accordingly so as to remain in $S^4$.  We denote the scaling factor
by $k_t^+(r)$.  It is easily computed.  Then we define a map $p_t^+$ of $S^4$
into $S^5$ by $$ p_t^+(v,r) = (k_t^+(r)v,(1-t)ir+t).  $$ In a similar way we
set $$ p_t^-(v,r) = (k_t^-(r)v,(1-t)ir-t).  $$ Then $t \mapsto \phi^+ \circ
p_t^+$ is a path of elements of ${\mathcal U}_3(C(S^4))$, as is $t \mapsto
\phi^- \circ p_t^-$.  We set $$ W_t = (\phi^+ \circ p_t^+)^*(\phi^- \circ
p_t^-)c_t $$ where $c_t = \left( \begin{array}{cc} I_2 & 0 \\ 0 & e^{\pi it}
\end{array} \right)$.  Then $W_t$ is a path of elements of ${\mathcal
U}_3(C(S^4))$ which goes from $U_2 \oplus I_1$ to $I_3$.

As indicated in lemma $8.5$ of \cite{Rie5}, a simple calculation shows that we
then obtain a ${\mathcal C}$-module isomorphism $\Phi$ from $X(U_0 \oplus I_1)$
to $X(I_3) = {\mathcal C}^3$ by setting $$ \Phi_F(t) = W_tW_0^{-1}F(t).  $$
Because $W_t$ is unitary, $\Phi$ preserves the ${\mathcal C}$-valued inner
products.  In particular, $\Phi$ will carry $X(U_0)$ and $X(U_0)^{\perp}$ to
orthogonal complementary submodules of ${\mathcal C}^3$.  Thus
$\Phi(X(U_0)^{\perp})$ will be a free rank-one submodule whose orthogonal
complement is not free, and so does not have a module basis.  Now in $X(U_0
\oplus I_1)$ the module $X(U_0)^{\perp}$ has as one choice of (normalized)
basis the evident constant function $$ E_3 = E_3(t) = \left( \begin{array}{c} 0
\\ 0 \\ 1 \end{array} \right) \in (C(S^4))^3.  $$ Thus $\Phi(X(U_0)^{\perp})$
will have as basis the $(C(S^4))^3$-valued function $$ N_0(t) = \Phi_{E_3}(t) =
W_tW_0^{-1}E_3(t) = W_tW_0^{-1}E_3, $$ which we can view as an element of
${\mathcal C}^3$.  (This element is closely related to our desired filter
$m_0$.)  But $W_0 = U_0 \oplus I_1$, and so $W_0^{-1}E_3 = E_3$.  Thus our
desired function is simply $$ N_0(t) = W_tE_3.  $$ Note that $N_0(1) = N_0(0)$
as needed, because $W_t$ goes from $U_2 \oplus I_1$ to $I_3$.  But if we look
at the definition of $W_t$, we see that we should calculate $$ ((\phi^- \circ
p_t^-)c_t)(v,r)E_3 = \phi^-(p_t^-(v,r))c_tE_3 $$ $$ = \left( \begin{array}{c}
k_t^-(r)v \\ e^{\pi it}((1-t)ir - t) \end{array} \right).  $$ Then to obtain
$N_0(t,v,r)$ we need only apply $(\phi^+ \circ p_t^+)^*$ to the above column
vector.  An explicit formula is not difficult to obtain, but we will not
display it here.  We note that $N_0(0,0,1) = (U_0(0,1) \oplus I_1)e_3 = e_3$,
where $e_3$ is the standard third basis vector of ${\mathbb C}^3$.  This is
related to the low-pass condition for our desired filter $m_0$.

Finally, we must relate all of the above to ${\mathbb T}^5$.  We view ${\mathbb
T}^5$ as ${\mathbb T} \times {\mathbb T}^4$, and we use the pinching map, $P$,
from ${\mathbb T}^4$ onto $S^4$.  To define it in a convenient way, we
momentarily let ${\mathbb T}$ be the interval $I = [-1,1]$ with the ends
identified.  Then $I^4$ is just the unit ball in ${\mathbb R}^4$ for the usual
supremum norm $\|\cdot\|_{\infty}$, and ${\mathbb T}^4$ is $I^4$ with certain
points of the boundary identified.  The pinching map $P$ then identifies all
the points of the boundary to just one point, forming $S^4$.  A convenient
formula for $P$, taking the boundary to $-e_5$, is $$ P(v) = (v\|v\|_2^{-1}
\sin(\pi\|v\|_{\infty}), \cos(\pi\|v\|_{\infty})) $$ for $v \in I^4$, where
$\|\cdot\|_2$ is the Euclidean norm on ${\mathbb R}^4$.  The first component of
the right-hand side must be viewed as an element of ${\mathbb C}^2$ under the
natural identification of ${\mathbb R}^4$ with ${\mathbb C}^2$.

We now let ${\tilde P} = I \times P$ to obtain a map from ${\mathbb T}^5$ onto
${\mathbb T} \times S^4$.  We use ${\tilde P}$ to pull back all of our earlier
construction to ${\mathbb T}^5$.  In particular we set $H_0 = N_0 \circ {\tilde
P}$.  Then $H_0 \in (C({\mathbb T}^5))^3 = {\mathcal A}^3$, and $\langle
H_0,H_0 \rangle_{\mathcal A} = 1$.  Furthermore, for $0 \in I^4$ we have $P(0)
= (0,1) \in S^4$, so that ${\tilde P}(0) = (0,0,1) \in {\mathbb T} \times S^4$.
We saw earlier that $N_0(0,0,1) = e_3$ and so $H_0(0) = e_3$ for $0 \in I^5$.
We need to have made our identification of ${\mathbb T}^5$ with $I^5$ in such a
way that $0$ goes to $0$.  Then $H_0$, as function on ${\mathbb T}^5$,
satisfies $H_0(0) = e_3$.

We must now compose $H_0$ with the ${\mathcal A}$-module isomorphism from
${\mathcal B} \cong C^*({\mathbb Z}^5)$ to ${\mathcal A}^3 \cong
(C^*(\Gamma))^3$ described early in this section.  But we must do so in such a
way that the third component in ${\mathcal A}^3$ corresponds to the submodule
of $C^*({\mathbb Z}^5)$ generated by the delta function at $0$ (discussed early
in this section), so that the function $h_0$ to which $H_0$ is carried
satisfies the low-pass condition $h_0(0) = 1$.  Of course it will also satisfy
$\langle h_0,h_0\rangle'_{\mathcal A} = 1$.  This is our desired filter.

We must check that this $h_0$ does not have two corresponding high-pass
filters.  If it did, then the (module-) orthogonal complement of the module
generated by $h_0$ would be free.  But this orthogonal complement corresponds,
and is isomorphic to, the submodule of ${\mathcal A}^3$ which is the pull back
of the module $X(U_0)$ over ${\mathcal C} = C({\mathbb T} \times S^4)$
discussed earlier.  Thus we need to know that this pulled-back module is not
free.  But it is not difficult to see that this module as essentially just the
module constructed from pulling $U_0$ itself back to $\mathbb T^4$.  That is,
we set $V_0 = U_0 \circ P$, so that $V_0 \in {\mathcal U}_2(C(\mathbb T^4))$.
Thus, just as argued above, 
we need to know that $V_0$ is not path connected to
a constant unitary, 
as this will imply that the $C(T^4)$-module $X(V_0)$ is not
free.  For the proof of this fact given below we have received much help from
Rob Kirby, from Elmer Rees through his reply to a query from Rob Kirby on our
behalf, from Jon Berrick, 
and especially from Jie Wu, who gave helpful answers
to many questions.

As is traditional, we let $[S^4,{\mathcal U}_2]$ denote the homotopy
equivalence classes of continuous functions from $S^4$ to ${\mathcal U}_2$
which preserve base points, and similary, for $[\mathbb T^4,{\mathcal U}_2],$
etc.  As discussed above, 
$U_0$ represents an element of $[S^4,{\mathcal U}_2]$
which is not trivial, in the sense that it is not homotopic to a constant
function.  By composing with the 
pinch map $P$ we obtain a mapping, $P^{\ast},$
from $[S^4,{\mathcal U}_2]$ to $[\mathbb T^4,{\mathcal U}_2].$ 
We need to check
that $P^{\ast}(U_0)$ is not the trivial element of $[\mathbb T^4,{\mathcal
U}_2].$ For this purpose it suffices to show that $P^{\ast}$ is injective.

Now ${\mathcal U}_2$ is a compact Lie group. It was shown by Milnor that any
nice topological group $G$ has a ``classifying space'', $Z,$ such that $G$ is
homotopic to $\Omega Z,$ where $\Omega Z$ is 
the loop-space of $Z$ (see Theorem
9.2.2 and the paragraph before Theorem 9.2.4 of \cite{sel}).  
We now use $Z$ to
denote any classifying space for ${\mathcal U}_2$ 
(the nature of the space will
not be of importance in the proof).  It is an elementary fact in homotopy
theory that forming loop spaces is the 
adjoint of forming (reduced) suspensions
(see Proposition 7.1.17 of \cite{sel}).  We denote suspensions by $\Sigma.$
Thus $$[S^4,{\mathcal U}_2]\;=\;[S^4,\Omega Z]\;=\;[\Sigma S^4, Z],$$ and a
similar equation holds with $S^4$ replaced by $\mathbb T^4.$ To $P^{\ast}$
there will then correspond the map $(\Sigma P)^{\ast},$ 
and our problem becomes
that of showing that $$(\Sigma P)^{\ast}:\;[\Sigma S^4,
Z]\;\rightarrow\;[\Sigma \mathbb T^4, Z]$$ is injective.  For this, it clearly
suffices to show that there is a mapping, $f,$ from $\Sigma S^4$ to $\Sigma
\mathbb T^4$ such that $\Sigma P\circ f$ is homotopic to the identity map on
$\Sigma S^4.$

We let $\wedge$ denote the standard smash product of spaces, as in Definition
6.2.12 of \cite{maun}, so that for any space $X$ we have by definition $\Sigma
X\;=\;\mathbb T\wedge X.$ For any two spaces $X$ and $Y$ we have that $\Sigma
(X\times Y)$ is homotopy equivalent to $\Sigma (X\wedge Y)\vee\Sigma
X\vee\Sigma Y,$ where $\vee$ denotes the reduced union, i.e.  one-point union,
or ``wedge" (see Exercise 15b of \cite{maun}).  When this exercise is applied
several times to $\Sigma \mathbb T^4$ we find that $\Sigma \mathbb T^4$ is
homotopy equivalent to a space of the form $\Sigma S^4\vee M$ for some space
$M$ which is a wedge of lower-dimensional spheres.  We let $f$ be the
composition of the evident map from $\Sigma S^4$ to $\Sigma S^4\vee M$ with a
homotopy equivalence from $\Sigma S^4\vee M$ to $\Sigma \mathbb T^4.$ We must
show that $\Sigma P\circ f$ is homotopic to the identity map on $\Sigma S^4.$

We first need some homology information about $P.$ Each of $\mathbb T^4$ and
$S^4$ is an orientable compact manifold, so $H_4(\mathbb
T^4)\;\cong\;\mathbb Z\;\cong\;H_4(S^4).$ The cube $I^4$ is homeomorphic with
the $4$-simplex, and the canonical map of $I^4$ onto $\mathbb T^4$ used above
is a generator for $H_4(\mathbb T^4).$ Similarly, the canonical map of $I^4$
onto $S^4$ which collapses the entire boundary is a generator of $H_4(S^4).$
But $P$ carries one of these canonical maps to the other.  Thus viewed as a map
in homology, $H_4(P)$ is an isomorphism from $H_4(\mathbb T^4)$ onto
$H_4(S^4).$ It follows from Theorem 4.4.10 of \cite{maun} that $$ H_5(\Sigma
P):\;H_5(\Sigma \mathbb T^4)\;\rightarrow\;H_5(\Sigma S^4) $$ is an
isomorphism.  Also, $H_5(\Sigma \mathbb T^4)\;\cong\; \mathbb
Z\;\cong\;H_5(\Sigma S^4)$ because $\Sigma S^4\;=\;S^5.$

Because $\Sigma \mathbb T^4$ is homotopy equivalent 
to $\Sigma S^4\vee M$ where
$M$ is a wedge of lower dimensional spheres, the map $f$ from $\Sigma S^4$ to
$\Sigma \mathbb T^4$ must give an isomorphism from $H_5(\Sigma S^4)$ to
$H_5(\Sigma \mathbb T^4).$ Hence $H_5(\Sigma P\circ f)$ is an isomorphism from
$H_5(\Sigma S^4)$ to itself.  
Because $\Sigma S^4\;=\;S^5,$ we have $H_n(\Sigma
S^4)=0$ except when $n=0$ or $n=5,$ where it is $\mathbb Z$ (see Theorem 4.6.6
of \cite{span}).  By the Hurewicz isomorphism theorem (Theorem 15.10 of
\cite{Ste}), $\pi_5(\Sigma S^4)\;=\;\pi_5(S^5)\;=\;H_5(S^5),$ with the
isomorphism being natural by the discussion in \cite{Ste} which preceeds 15.10.
Thus $\pi_5(\Sigma P\circ f)$ is an isomorphism.  It preserves orientation, so
$\Sigma P\circ f$ is homotopic to the identity map on $\Sigma S^4$.  This
implies, as described above, that $(\Sigma P)^{\ast}$ is injective, so that
$P^{\ast}$ is injective.  This in turn implies that $P^{\ast}(U_0)$ is
non-trivial in $[\mathbb T^4,{\mathcal U}_2],$ as we promised to show.

Finally, we remark that, much as discussed at the end of the proof of Theorem
3.1, $h_0$ can be approximated arbitrarily closely by low-pass filters which
are smooth (i.e., infinitely differentiable).  For sufficiently close
approximates they will be connected to $h_0$ by paths through the space of
low-pass filters.  This can be seen directly, or proposition $5.1$ of
\cite{Rie4} can be used.  Such smooth approximates will have complementary
modules isomorphic to that of $h_0$, and so again not free.  We suspect that in
a similar way our specific $h_0$ can be adjusted so that it satisfies condition
(iii) of Definition 2.1, but we have not worked out the details.

\end{document}